\documentclass[12pt,TimesNewRoman]{article} 
\usepackage{pgf,tikz}
\usepackage[utf8]{inputenc} 
\usepackage[T1]{fontenc}    
\usepackage{amsmath}
\usepackage[margin=1in,footskip=0.25in]{geometry}
\usepackage{multirow}	
\usepackage[utf8]{inputenc}	
\usepackage{amssymb}	
\usepackage{amsfonts}
\usepackage{latexsym}
\usepackage{tipa}
\usepackage{verbatim}
\usepackage{amsthm}
\usepackage{mathrsfs}
\usepackage[all]{xy}
\usepackage{graphics}
\usepackage{multicol}
\usepackage{xcolor,rotate}

\def\->{{{\rightarrow}}}

\newtheorem{theorem}{Theorem}

\theoremstyle{definition}
\newtheorem{definition}[theorem]{Definition}
\newtheorem{example}[theorem]{Example}

\newtheorem{re}{Remark}

\date{}
\title{Generalized Fractional Differential Ring}

\author{Zeinab Toghani\:, Luis Gaggero \\
Queen Mary University of London \\
Universidad Autónoma del Estado de Morelos(CIICAp)}
 
\begin{document}
\maketitle
\let\thefootnote\relax\footnotetext{
\emph{Emails: toghaniii@gmail.com, lgaggero@uaem.mx}. Corresponding Author: Luis Gaggero-Sager.}
\begin{abstract}
 There are many possible definitions of derivatives, here we present some and present one that we have called generalized that allows us to put some of the others as a particular case of this
  but, what interests us is to determine that there is an infinite number of possible definitions of fractional derivatives, all are correct as differential operators each of them must be properly defined its algebra. 
 
We introduce a generalized version of fractional derivative that
extends the existing ones in the literature. To those extensions it is
associated a differentiable operator and a differential ring and applications that shows the advantages of the
generalization. 

We also review the different definitions of fractional derivatives
proposed by  Michele Caputo in \cite{GJI:GJI529},
Khalil, Al Horani, Yousef, Sababheh in \cite{khalil2014new}, Anderson and Ulness
in \cite{anderson2015newly},
Guebbai and Ghiat in [GG16], Udita N. Katugampola in \cite{katugampola2014new}, Camrud in \cite{camrud2016conformable}
and it is shown how the generalized version contains the previous ones
as a  particular cases.
\end{abstract}


\section{Introduction}
Fractional derivative was defined for responding to a question' what does it mean $ \frac{d^{\alpha}f}{dt^{\alpha}} $ if $ \alpha=\frac{1}{2} $ ' in 1695. 
  Following that, finding the right  definition of fractional derivative has
attracted significant attention of researcher and in the last few years it
has seen significantly progress in mathematical and non-mathematical
journals(see \cite{laskin2002fractional}, \cite{wei2015some}, \cite{agarwal2012applications}, \cite{du2013measuring}, \cite{cenesiz2015solutions}, \cite{wang2017evaluation}, \cite{abu2014fractional}, \cite{namias1980fractional}, \cite{younis2013travelling}, \cite{anderson2016taylor}). In fact, there are articles which in few months have gained hundreds of citations. In particular in past three years several definitions of fractional derivative have been proposed(see \cite{khalil2014new}, \cite{anderson2015properties}, \cite{sousa2017m}, \cite{camrud2016conformable}, \cite{miller1993introduction}, \cite{Ortigueira:2015:FD:2777696.2777939}, \cite{katugampola2014new}, \cite{iyiola2016some}, \cite{abdeljawad2015conformable}, \cite{abdeljawad2015conformable}, \cite{guebbai2016new}, \cite{anderson2016taylor}, \cite{sarikaya2017new}, \cite{caputo2015new}, \cite{du2013measuring}, \cite{agarwal2012applications}, \cite{laskin2000fractional}). Since some of previous definitions do not satisfy the classical
formulas of the usual derivative, it has been proposed an ad hoc
algebra associated to each definition. To unify that diversity, we
propose a version of fractional derivative that has the advantages
that generalized the already existing in the literature and where the
different algebras are unified under the notion of fractional
differential ring.
  
  The present paper is organized as
follows: In the second section we give the previous definitions of fractional derivative and our generalized fractional derivative(GFD) definition, in the third section we introduce a fractional differential ring, in the forth section we give some result of GFD, in the fifth section we give a definition of partial fractional differential derivative, in the sixth section we give a definition of GFD when $ \alpha\in(n,n+1] $. 
\section{Fractional derivative}
Let $ \alpha\in (0,1] $ be a fractional number, we want to give a definition of generalized fractional derivative of order $ \alpha$ for a differentiable function $ f $.
We denote $ \alpha-$th derivative of $ f$ by $ D^{\alpha}(f) $ and we denote the first derivative of $ f$ by $ D(f) $. 

We begin the present section listing previous definition of fractional
derivative; later we present our proposal of generalized one showing
how it contains the once already described. We finish the section
providing some examples.
\begin{enumerate}
\item The Caputo fractional derivative was defined by Michele Caputo in \cite{caputo1969elasticita}:
\begin{equation}\label{caputo}
 D^{\alpha}(f)=\frac{1}{\Gamma(1-\alpha)}\int_{a}^{t}\frac{f'(x)}{(t-x)^{\alpha}}dx.
\end{equation}
\item The conformable fractional derivative was defined by Khalil, Al Horani, Yousef and Sababheh in \cite{khalil2014new} : 
\begin{equation}\label{arabe}
 D^{\alpha}f(t)=\lim_{\varepsilon\rightarrow 0}\frac{f(t+ t^{1-\alpha}\varepsilon)-f(t)}{\varepsilon}.
\end{equation}
 \item The conformable fractional derivative was defined by Anderson and Ulness in \cite{anderson2015newly}: 
\begin{equation}\label{gringos}
 D^{\alpha}f(t)=(1-\alpha)\mid t \mid^{\alpha}f(t)+\alpha \mid t\mid^{1-\alpha}Df.
 \end{equation}
 \item The fractional derivative was defined by Udita N.Katugampola in \cite{katugampola2014new}:
 \begin{equation}\label{algeriano}
 D^{\alpha}f(t)=\lim_{\varepsilon\rightarrow 0}\frac{f(te^{\varepsilon t^{-\alpha}})-f(t)}{\varepsilon}.
 \end{equation}
\item The fractional derivative was defined by Guebbai and Ghiat in \cite{guebbai2016new} for an increasing and positive function $ f $:
\begin{equation}\label{algeriano}
 D^{\alpha}f(t)=\lim_{\varepsilon\rightarrow 0}\left(\frac{f(t+f(t)^{\frac{1-\alpha}{\alpha}}\varepsilon)-f(t)}{\varepsilon}\right)^{\alpha}.
 \end{equation}
 \item The conformable ratio derivative was defined by Camrud in \cite{camrud2016conformable} for a function $f (t)\geq 0$ with $Df(t)\geq 0$:
 \begin{equation}
 D^{\alpha}f(t)=\lim_{\varepsilon\rightarrow 0}f(t)^{1-\alpha}\left(\frac{f(t+\varepsilon)-f(t)}{\varepsilon}\right)^{\alpha}.
 \end{equation}
 \end{enumerate}
From all these definitions, we propose a definition that unifies almost all of them.
 \begin{definition}\label{tt}
Given a differentiable function $f :[0,\infty) \rightarrow \mathbb{R}$, \textbf{the generalized fractional derivative}(GFD) for $ \alpha \in (0,1] $ at point $ t $ is defined by :
\begin{displaymath}
\begin{array}{llll}
 D^{\alpha}f(t)=
 \lim_{\varepsilon\rightarrow 0}\frac{f(t+w_{t,\alpha} t^{1-\alpha}\varepsilon)-f(t)}{\varepsilon},
 \end{array}
 \end{displaymath}
 where $w_{t,\alpha}$ is a function that may depend on $ \alpha $ and $ t $.
\end{definition}

\begin{re}\label{AA}
As a consequence of definition \ref{tt} we can see 
$$ D^{\alpha}f(t)=w_{t,\alpha} t^{1-\alpha}Df(t). $$
\end{re}

\begin{definition}
A differentiable function $ f: [0,\infty) \rightarrow \mathbb{R}$ is said to be \textbf{$\alpha-$generalized fractional differentiable function} over $ [0,\infty)  $ if it exists $ D^{\alpha}(f)(t) $ for all $ t\in [0,\infty) $ for $ \alpha\in (0,1] $.
\end{definition}
We denote $ C^{\alpha}[0,\infty) $ the set of $ \alpha- $generalized differentiable functions with real values in the interval $[0,\infty)$ in variable $t$. The set $ (C^{\alpha}[0,\infty),+,.)$ is a ring.
In the following we want to see the relation between GFD and the others definitions:
\begin{enumerate}
\item The fractional derivative of Khalil, Al Horani, Yousef and Sababheh in \cite{khalil2014new} is a particular case of GFD  where $ w_{t,\alpha}=1$.  
\item The fractional derivative of Anderson and Ulness in \cite{anderson2015newly} is a particular case of GFD where 
\[
w_{t,\alpha}=\frac{(1-\alpha)t^{\alpha}f(t)+\alpha t^{1-\alpha}Df}{\alpha t^{1-\alpha}}. 
\]
In this fractional derivative $ w_{t,\alpha} $ depends on $ \alpha $ and $t$.
\item The fractional derivatives of Guebbai and Ghiat in  \cite{guebbai2016new} and Camrud in  \cite{camrud2016conformable} are particular cases of GFD where $ w_{t,\alpha}=\left(\frac{tDf}{f} \right)^{\alpha-1}$.
\end{enumerate}
 We are particularly interested in discussing GFD where $w_{t,\alpha}=g(t,\alpha)\tau^{\alpha-1}$ such that $  g:[0,\infty)\times(0,1]\rightarrow \mathbb{R}$ is
 a function and $ \tau$ is the characteristic of system 
with the properties
\[
w_{t,\alpha} =1 \quad \textit{if and only if} \quad \alpha=1.
\]
If the system is periodic with period $ T$, then we have $ \tau=T $. In the quantum systems $ \tau $ is the Bohr radius and in astronomy $ \tau $ is the light year. The characteristic of system $ \tau $ depends on the systems and the derivative. If $ t $ is time, $ \tau  $ is time too. If $ t $ is space, $ \tau  $ is space too. In fact the unit of $ t $ is $ \tau $, i.e., $ t=c\tau $ where $ c $ is a constant. In the general $ \tau =1$.
\begin{example}
Let $ \alpha,\beta \in (0,1] $. Let $ f,h $ be two functions in $C^{\alpha}[0,\infty)  $. We suppose $w_{t,\alpha}= g(t,\alpha)\tau^{\alpha-1} $ with $ \tau=1$.
\begin{enumerate}
 \item  If $g(t,\alpha)$ is a function with $ g(0,0)=0,$ then $ \lim_{\alpha \rightarrow 0} D^{\alpha}(f)=0$. 
\item If $ g(t,\alpha)=\alpha $, we have the chain rule
$$ D^{\alpha}(f\circ h) =\frac{t^{\alpha-1}}{\alpha} D^{\alpha}\left( f(h)\right) D^{\alpha}(h).$$
\end{enumerate} 
 \end{example}
 
\begin{example}
We want to present the corresponding figure to the generalized fractional derivatives for $ \alpha=\frac{3}{4}$ for a trigonometric, using all the fractional derivative definitions that we have already mentioned in this article. It can be seen from all the figures that in principle these definitions do not find a reason to discard them. That is, they have a fairly reasonable behavior. We consider $ f(t)=\sin(2t) $,
the figure of $f$ can be seen in the following picture.
 \begin{figure}[h]
\centering\label{expo}
\includegraphics[height=5cm]{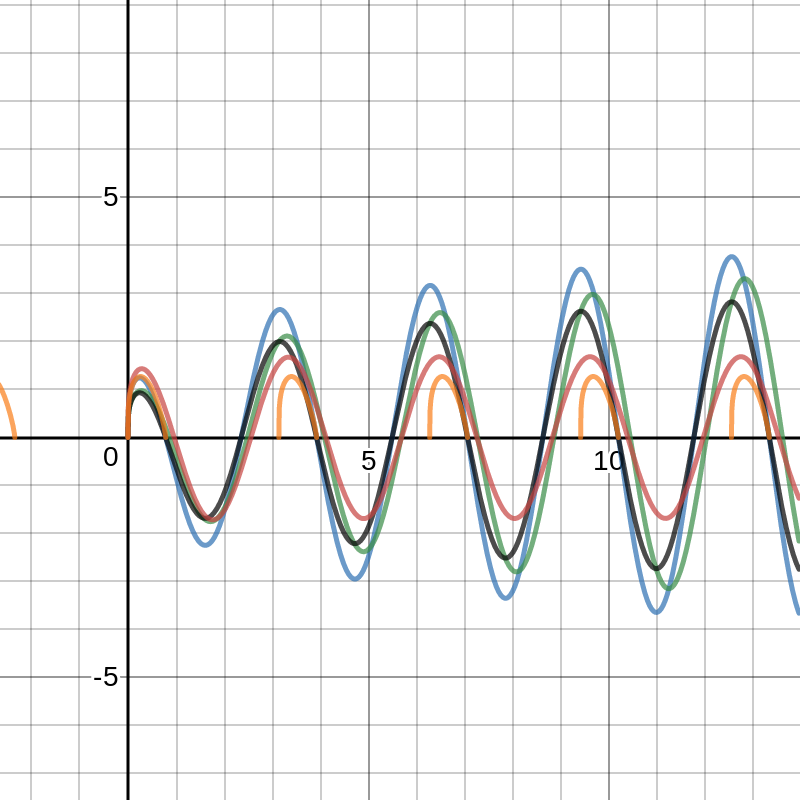}
\caption{red:Caputo, green:Khalil et al, blue:Anderson et al , orange:Guebbai et al, black:GFD when $ w_{t,\alpha}= \alpha$.}
\end{figure}
\end{example}

\section{Generalized fractional differential ring}
In this section we want to stress out that instead of defining a new derivative, we focus on the notion of differentiable operator and the ring that it carries with.
\begin{definition}
Let $R$ be a commutative ring with unity. A derivation on $R$ is a map $d : R \rightarrow R$ that satisfies
$d(a + b) = d(a) + d(b)$ and $d(ab) = d(a)b + ad(b),\forall a, b \in R$. The pair $(R, d)$ is called a differential ring(see \cite{Ritt:1950}).
\end{definition}
  
 \begin{theorem}
 Let $ \alpha\in (0,1] $. The ring $ C^{\alpha}[0,\infty)$ with operator $D^{\alpha}: C^{\alpha}[0,\infty)\rightarrow C^{\alpha}[0,\infty)$ is a differential ring.  
\end{theorem}
 \begin{proof}
Since $ C^{\alpha}[0,\infty)$ is a commutative ring with unity $ f(t)=1 $ and the derivation $D^{\alpha}$ for $ \alpha\in (0,1] $ satisfies following properties from remark \eqref{AA}
 \begin{enumerate}
\item $ D^{\alpha}(af_{1}+bf_{2})=aD^{\alpha}(f_{1})+bD^{\alpha}(f_{2}),  $  $ \forall f_{1},f_{2} \in C^{\alpha}[0,\infty), \forall a,b\in \mathbb{R}$,
\item$ D^{\alpha}(f_{1}f_{2})=f_{1}D^{\alpha}(f_{2})+f_{2}D^{\alpha}(f_{1})$,  $ \forall f_{1},f_{2} \in C^{\alpha}[0,\infty) $. 
  \end{enumerate} 
  \end{proof}
Let $ \alpha\in (0,1]$ be a fractional number, let $ f_{1},f_{2} \in C^{\alpha}[0,\infty)$ be two functions, GFD has the following properties:
\begin{enumerate}
 \item 
 $D^{\alpha}(\frac{f_{1}}{f_{2}})=\frac{f_{2}D^{\alpha}f_{1}-f_{1}D^{\alpha}f_{2}}{f_{2}^2}.$
 \item $
D^{\alpha}(f_{1}\circ f_{2}) =\frac{t^{\alpha-1}}{w_{t,\alpha}} D^{\alpha}\left( f_{1}(f_{2})\right)D^{\alpha}(f_{2}).$
\item  $D^{\alpha+\beta}(f_{1})=\frac{w_{t,\alpha}w_{t,\beta}t}{w_{t,\alpha+\beta}}D^{\alpha}D^{\beta}(f_{1}).$
 \end{enumerate}
It is easy to see these properties from remark \eqref{AA}. If $ w_{t,\alpha}=t^{1-\alpha} $ we have the equality
$$ D^{\alpha}(f_{1}\circ f_{2}) =D^{\alpha}\left( f_{1}(f_{2})\right) D^{\alpha}(f_{2}).$$
If $ \frac{w_{t,\alpha+\beta}}{w_{t,\alpha}w_{t,\beta}}=t$ we have the equality
$$ D^{\alpha+\beta}(f_{1})= D^{\alpha}D^{\beta}(f_{1}), \forall  \alpha,\beta\in (0,1].$$
Parts 4 and 5 of the properties imply that we can create function spaces with different algebras using different expressions for $w_{t,\alpha}$.

By considering previous properties of GFD we called $ C^{\alpha}[0,\infty)$ a \textbf{$w_{t,\alpha} -$generalized fractional differential ring} of functions and we denote it by $(C^{\alpha}[0,\infty) ,D^{\alpha},w_{t,\alpha}) $.
Let $ I \subset  C^{\alpha}[0,\infty)$ be an ideal. If $ D^{\alpha}(I)\subset I $ then the ideal $ I $ is called a \textbf{$w_{t,\alpha} -$generalized fractional differential ideal}.
By using the previous properties we can see the following result:
  \begin{theorem}
 Let $ \alpha\in(0,1] $. Associated to any $ \alpha $ and any $ w_{t,\alpha} $ there exists a fractional
differential ring.
 \end{theorem}
 \section{Some results of Generalized Fractional Derivative} 
Let $ \alpha\in (0,1] $ be a fractional number and $ t\in [0,\infty) $ then GFD has the following properties:
 \begin{enumerate}
 \item
$ D^{\alpha}(\frac{t^{\alpha}}{\alpha w_{t,\alpha}})=1 $,
\item
$ D^{\alpha}(\sin(\frac{t^{\alpha}}{\alpha w_{t,\alpha}}))=\cos(\frac{t^{\alpha}}{\alpha w_{t,\alpha}})$,
\item
$ D^{\alpha}(\cos(\frac{t^{\alpha}}{\alpha w_{t,\alpha}}))=-\sin(\frac{t^{\alpha}}{\alpha w_{t,\alpha}})$,
\item
$ D^{\alpha}(e^{(\frac{t^{\alpha}}{\alpha w_{t,\alpha}})})=e^{(\frac{t^{\alpha}}{\alpha w_{t,\alpha}})}$.
 \end{enumerate}
\begin{theorem}(Rolle's theorem for $ \alpha- $Generalized Fractional Differentiable Functions)

Let $ a>0  $, let $ f:[a,b]\rightarrow \mathbb{R} $ be a function with the properties that
\begin{enumerate}
\item $ f $ 
is continuous on $ [a,b] $,
\item $f$ is $ \alpha- $generalized fractional differentiable on $(a, b)$ for some $ \alpha\in(0,1] $, 
\item  $ f(a)=f(b) . $
\end{enumerate}
Then, there exist $ c\in (a,b) $ such that $ D^{\alpha}f(c)=0 $.
\end{theorem}
\begin{proof}
Since $f  $ is continuous on $ [a,b]$ and $ f(a)=f(b)  $, then the funtion $ f $ has a local extreme in a point $ c\in (a,b) $. Then  
$$ D^{\alpha} f(c)=\lim_{\varepsilon\rightarrow 0^{+}}\frac{f(c+w_{t,\alpha} c^{1-\alpha}\varepsilon)-f(c)}{\varepsilon}=\lim_{\varepsilon\rightarrow 0^{-}}\frac{f(c+w_{t,\alpha} c^{1-\alpha}\varepsilon)-f(c)}{\varepsilon}.$$ 
But two limits have different signs. Then $D^{\alpha}f(c)=0 .  $ 
\end{proof}
\begin{theorem}(Mean Value Theorem for $ \alpha- $Generalized Fractional Differentiable Functions)
Let $ a>0 $ and $ f:[a,b]\rightarrow \mathbb{R} $ be a function with the properties that
\begin{enumerate}
\item $ f $ 
is continuous on $ [a,b] $,
\item $f$ is $ \alpha-$Generalized fractional differentiable on $(a,b)$ for some $ \alpha\in (0,1] $.
\end{enumerate}
Then, there exists $ c\in (a,b) $ such that $ D^{\alpha}f(c)=\frac{\alpha w_{t,\alpha}(f(b)-f(a))}{b-a}$.
\end{theorem}
\begin{proof}
Consider function  
\[
h(t)=f(t)-f(a)-\frac{\alpha w_{t,\alpha}(f(b)-f(a))}{b-a}\left(\frac{t^{\alpha}}{\alpha w_{t,\alpha}} -\frac{a^{\alpha}}{\alpha w_{t,\alpha}}\right) ,
\]
Then, the function $ h $ satisfies the conditions of the fractional Rolle’s theorem. Hence, there exists $ c\in (a,b) $ such that $ D^{\alpha}h(c)=0. $ We have the result since 
\[
 D^{\alpha}h(c)= D^{\alpha}f(c)-\frac{\alpha w_{t,\alpha}(f(b)-f(a))}{b-a}(1)=0
\] 
\end{proof}
%
%
 \section{Generalized Partial Fractional Derivative}
In this section we introduce a partial fractional derivative of first and second order, also we introduce a partial fractional differential ring. 
\begin{definition}\label{tttt}
Let $f(t_{1},\cdots,t_{n}):[0,\infty)^n\rightarrow \mathbb{R}$ be a function with $ n  $ variables such that $\forall i$, there exists the partial derivative of $ f $ respect to $t_{i}$. Let $ \alpha\in (0,1]$ be a fractional number. We define $ \alpha-$\textbf{generalized partial fractional derivative}(GPFD) of $ f $ with respect to $ t_{i} $ at point $ t=(t_{1},\ldots,t_{n}) $
\begin{displaymath}
\begin{array}{llll}
 \frac{\partial^{\alpha}f(t)}{\partial t_{i}^{\alpha}}=
 \lim_{\varepsilon\rightarrow 0}\frac{f(t_{1},\ldots,t_{i}+w_{t_{i},\alpha} t_{i}^{1-\alpha}\varepsilon,\ldots,t_{n})-f(t)}{\varepsilon},
 \end{array}
 \end{displaymath}
 where $w_{t_{i},\alpha}$ can be a function depend on $ \alpha$ and $ t_{i} $.
 \end{definition}
\begin{re}\label{BB}
As a consequence of definition \eqref{BB} we can see for $ \alpha\in(0,1] $ and $ 1\leq i\leq n $:
\[
 \frac{\partial^{\alpha}f}{\partial t_{i}^{\alpha}}(t)=w_{t_{i},\alpha}( t_{i})^{1-\alpha} \frac{\partial f}{\partial t_{i}}(t).
\]
\end{re}
Let $ \alpha\in (0,1] $ and $ 1\leq i\leq n $. A partial differentiable function $ f: [0,\infty)^n \rightarrow \mathbb{R}$ is said to be a $ \alpha- $\textbf{generalized partial fractional differentiable function} respect to $ t_{i} $ over $ [0,\infty)  $ if exists $ \frac{\partial^{\alpha}f(t)}{\partial t_{i}^{\alpha}} $ for all $ t\in [0,\infty) $. We denote by $ C_{i}^{\alpha}[0,\infty)^n $ the set of $ \alpha-$generalized partial fractional differentiable functions respect to $ t_{i} $ with real values in the interval $[0,\infty)^n$ in variable $t=(t_{1},\ldots,t_{n})$. The set $ (C_{i}^{\alpha}[0,\infty)^n,+,.)$ is a ring.
  \begin{theorem}
Let $ \alpha\in(0,1] $ and $ 1\leq i\leq n $. The ring 
$ C_{i}^{\alpha}[0,\infty)^n$ with operator $$ \frac{\partial^{\alpha}}{\partial t_{i}^{\alpha}}: C_{i}^{\alpha}[0,\infty)^n\rightarrow  C_{i}^{\alpha}[0,\infty)^n$$
 is a differential ring.
\end{theorem}
\begin{proof}
Since the ring $ C_{i}^{\alpha}[0,\infty)^n$ is a commutative ring with unity $ f(t_{1},\ldots,t_{n})=1 $ and the derivation $\frac{\partial^{\alpha}}{\partial t_{i}^{\alpha}}$ for $\alpha\in (0,1] $ satisfies the following properties from remark \eqref{BB}; 
\begin{enumerate}
\item
 $\frac{\partial^{\alpha}(f_{1}+f_{2})}{\partial t_{i}^{\alpha}} =\frac{\partial^{\alpha}(f_{1})}{\partial t_{i}^{\alpha}}+\frac{\partial^{\alpha}(f_{2})}{\partial t_{i}^{\alpha}}\quad f_{1},f_{2}\in C^{\alpha}[0,\infty)^n,$
\item $\frac{\partial^{\alpha}(f_{1}f_{2})}{\partial t_{i}^{\alpha}} =f_{1}\frac{\partial^{\alpha}f_{2}}{\partial t_{i}^{\alpha}}+f_{2}\frac{\partial^{\alpha}f_{1}}{\partial t_{i}^{\alpha}}\quad f_{1},f_{2}\in C^{\alpha}[0,\infty)^n$.
 \end{enumerate}
\end{proof}
Let $ \alpha\in (0,1]$ and $ 1\leq i\leq n $, let $ f_{1},f_{2} \in C^{\alpha}[0,\infty)^n$ be two functions, then GPFD has the following properties from remark \eqref{BB}:
\begin{enumerate}
 \item 
$\frac{\partial^{\alpha}(\frac{f_{1}}{f_{2}})}{\partial t_{i}^{\alpha}}=\frac{f_{2} \frac{\partial^{\alpha}f_{1}}{\partial t_{i}^{\alpha}}-f_{1} \frac{\partial^{\alpha}f_{2}}{\partial t_{i}^{\alpha}}}{f_{2}^2},$
\item $
 \frac{\partial^{\alpha}f_{1}\circ f_{2}}{\partial t_{i}^{\alpha}}= \frac{t_{i}^{\alpha-1}}{w_{t_{i},\alpha}} 
  \frac{\partial^{\alpha}( f_{1}(f_{2}))}{\partial t_{i}^{\alpha}}\frac{\partial^{\alpha}(f_{2})}{\partial t_{i}^{\alpha}},
 $
\item 
$ \frac{\partial^{\alpha+\beta}( f_{1})}{\partial t_{i}^{\alpha}}=\frac{w_{t_{i},\alpha}w_{t_{i},\beta}t_{i}}{w_{t_{i},\alpha+\beta}} \frac{\partial^{\alpha}}{\partial t_{i}^{\alpha}}\frac{\partial^{\alpha}(f_{1})}{\partial t_{i}^{\alpha}}.$
 \end{enumerate}
By considering previous properties of GPFD we called the ring $ C_{i}^{\alpha}[0,\infty)^n$ a $ w_{t_{i},\alpha}- $ \textbf{generalized partial fractional differential ring}. We denote it by $ (C_{i}^{\alpha}[0,\infty)^n,\frac{\partial^{\alpha}}{\partial t_{i}^{\alpha}} ,w_{t_{i},\alpha})  $.
 
 We can see the following result by using the previous properties.
 \begin{theorem}
 Let $ \alpha\in(0,1] $ and $ 1\leq i \leq n $. Associated to any $ \alpha $ and any $ w_{t_{i},\alpha} $ there is a partial fractional
differential ring.
\end{theorem}
 \begin{example}
Let $ f(t_{1},t_{2}) =t_{1}^3\sin(t_{2})$, let $ \alpha\in [0,1)$. We have 
$$ \frac{\partial^{\alpha}f}{\partial t_{1}^{\alpha}}= w_{t,\alpha}(t_{1})^{1-\alpha}(3t_{1}^2)\sin(t_{2}).$$
\end{example}
\begin{definition}\label{dd}
Let $ \alpha\in(0,1]$ be a fractional number. We define $ \alpha- $\textbf{generalized partial fractional derivative of second order} with respect to $t_{i}$ and $ t_{j} $ at point $ t=(t_{1},\cdots,t_{n}) $ is 
\begin{displaymath}
\begin{array}{llll}
  \frac{\partial^{\alpha^2}f(t)}{\partial t_{j}^{\alpha} \partial t_{i}^{\alpha}}=\frac{\partial^{\alpha}}{\partial t_{j}^{\alpha}}(\frac{\partial^{\alpha}f(t)}{\partial t_{i}^{\alpha}})=
 \lim_{\varepsilon\rightarrow 0}\frac{ \frac{\partial^{\alpha}f}{\partial t_{i}^{\alpha}}(t_{1},\ldots,t_{j}+w_{t_{j},\alpha} t_{j}^{1-\alpha}\varepsilon,\ldots,t_{n})-\frac{\partial^{\alpha}f(t)}{\partial t_{i}^{\alpha}}}{\varepsilon}.
 \end{array}
 \end{displaymath}
\end{definition}
 \begin{re}\label{CC}
As a consequence of definition \eqref{dd} we can see for $ \alpha\in(0,1] $ and $ 1\leq i,j\leq n $;
\[
  \frac{\partial^{\alpha^2}f(t)}{\partial t_{j}^{\alpha} \partial t_{i}^{\alpha}}=w_{t_{j},\alpha}w_{t_{i},\alpha}(t_{j} t_{i})^{1-\alpha} \frac{\partial }{\partial t_{j}}(  \frac{\partial f}{\partial t_{i}}(t)).
\]
\end{re}
A partial differentiable function of second order $ f: [0,\infty)^n \rightarrow \mathbb{R}$ is said to be a $ \alpha- $generalized fractional partial differentiable function of second order respect to $ t_{i}$ and $t_{j} $ over $ [0,\infty)  $ if $  $ exists $ \frac{\partial^{\alpha^2}f(t)}{\partial t_{j}^{\alpha}\partial t_{i}^{\alpha}}$ for all $ t\in [0,\infty) $. We denote $ C_{i,j}^{\alpha^2}[0,\infty)^n $ the set of $ \alpha-$generalized partial fractional differentiable functions of second order respect to $ t_{i}$ and $t_{j} $ with real values in the interval $[0,\infty)^n$ in variable $t=(t_{1},\ldots,t_{n})$. The set $ (C_{i,j}^{\alpha^2}[0,\infty)^n,+,.)$ is a ring.
\section{Generalized Fractional Derivative for $ \alpha\in (n,n+1] $}
In this section we define a fraction differential derivative for $ \alpha\in (n,n+1]. $ 
\begin{definition}\label{ttt}
Let $ \alpha\in (n,n+1]$ be a fractional number for $ n\in \mathbb{N} $, let $ f:[0.\infty) \rightarrow \mathbb{R}$ be a $ n -$ differentiable. The \textbf{generalized fractional derivative} of order $ \alpha $ is defined by 
\begin{displaymath}
\begin{array}{llll}
 D^{\alpha}f(t)=
 \lim_{\varepsilon\rightarrow 0}\frac{f^{[\alpha]-1}(t+w_{t,\alpha}  t^{[\alpha]-\alpha}\varepsilon)-f^{[\alpha]-1}(t)}{\varepsilon}.
\end{array}
 \end{displaymath}
where $[\alpha]$ is the smallest integer greater than or equal to $ \alpha. $ 
 \end{definition}
 As a consequence of definition \ref{ttt} we can see 
\[
D^{\alpha}(f)=w_{t,\alpha}t^{[\alpha]-\alpha}D^{[\alpha]}(f),
\]
where $ \alpha\in (n,n+1] $.

Let $ n<\alpha\leq n+1 $. A function $ f: [0,\infty) \rightarrow \mathbb{R}$ is said to be $\alpha-$generalized differentiable over $ [0,\infty)$ if there exists $ D^{\alpha}(f)(t) $ for all $ t\in [0,\infty). $ We denote $ C^{\alpha}[0,\infty) $ the set of $ \alpha- $generalized fractional differentiable functions with real values in the interval $[0,\infty)$ in variable $t$. The set $ (C^{\alpha}[0,\infty),+,.)$ is a ring.
 \begin{theorem}
The ring $ (C^{\alpha}[0,\infty),+,.)$ with operator $ D^{\alpha} $ is not a differential ring for fractional number $ \alpha\in (n,n+1] $ .
\end{theorem}
\begin{proof}
Since $ D^{\alpha}(fg)\neq fD^{\alpha}g+gD^{\alpha}(f) $ for every $ f,g\in C^{\alpha}[0,\infty) $.
\end{proof}
\section{$\alpha$-Fractional Taylor Series}
There are some articles about fractional Taylor series see( \cite{anderson2015properties} , \cite{usero2008fractional},\cite{munkhammar2004riemann}, \cite{yang2011generalized}).
In this section we use GFD to define a fractional taylor series for a function $ f\in C^{r}[0,\infty)$ for every fractional number $ r $.

Let $0 <\alpha<1 $, we define the $\alpha$-fractional taylor series of $ f $ at real number $ x_{0} $ 
\[
f(x)=f(x_{0})+\sum_{i=1}^{\infty} \frac{D^{i}f(x_{0})}{w_{x,\alpha}\overline{(\alpha+i-1)!}}  (x-x_{0})^{\alpha+i-1} ,
\]
where $\overline{(\alpha+i-1)!}= \alpha(\alpha+1)\cdots(\alpha+i-1) $.

Let $ 1<\alpha\leq 2 $, we define the $\alpha$-fractional taylor series of $ f $ at real number $ x_{0}$ 
\[
f(x)=f(x_{0})+Df(x_{0})(x-x_{0})+\sum_{i=2}^{\infty}\frac{D^{i}f(x_{0})(x-x_{0})^{\alpha+i-2}}{w_{x,\alpha}\overline{(\alpha+i-2)!}},
\]
where $\overline{(\alpha+i-2)!}= (\alpha-1)\alpha(\alpha+1)\cdots(\alpha+i-2) .$
%
 Let $ n<\alpha\leq n+1 $ such that $ \alpha=n+A $ with $ 0<A<1 $ we define the $\alpha$-fractional taylor series of $ f $ at real number $ x_{0}$, 
 \[
f(x)=f(x_{0})+\sum_{i=1}^{n}\frac{D^if(x_{0})(x-x_{0})^{i}}{i!}+\sum_{i=n+1}^{\infty}\frac{D^{i}f(x_{0})}{w_{x,\alpha}\overline{(A+i-1)!}}(x-x_{0})^{A+i-1},
\]
where $ \alpha=n+A $, $\overline{(A+i-1)!}= A(A+1)\cdots(A+i-1) $.
\section{Application to differential equations}
There are some articles for applications of fractional differential derivative such as \cite{younis2013travelling}, \cite{agarwal2012applications} ,\cite{abu2014fractional}, \cite{wang2017evaluation}.
In this section we solve some (partial)fractional differential equations by using our definitions.
At the first we solve the fractional differential equations with the form 
\begin{equation}\label{SS}
 aD^{\alpha}y+by=c,
 \end{equation} 
where $y=f(t)$ be a differentiable function and $0<\alpha<1  $. 

By substituting GFD in the equation \eqref{SS} we have
\[
aw_{t,\alpha}t^{1-\alpha}Dy+by=c\Longrightarrow Dy+\frac{bt^{\alpha-1}}{aw_{t,\alpha}}y=\frac{t^{\alpha-1}c}{aw_{t,\alpha}},  
\]
the solutions of this equation have the form $y(t)= \frac{c}{b}+c_{1}e^{(-\frac{bt^{\alpha}}{aw_{t,\alpha}\alpha})}. $
\begin{example}
We consider the partial fractional differential equation with boundary conditions
\begin{equation}\label{tropsol2}
\begin{cases}
u_{t}+2\sqrt[3]{xu_{x}}+u=x^{2},t>0\\
u(t,0)=0,\\
u(0,x)=0,
\end{cases}
\end{equation}
where $ u(x,t)$ be a differentiable function respect to $ x $ and $ t $,  $ u(x,t)$ be a   
 $ \frac{1}{3}- $partial fractional differentiable function of first order respect to $x $ 
, $u_{t}=\frac{\partial u}{\partial t } $ and $\sqrt[3]{u_{x}}=\frac{\partial^{\frac{1}{3}} u}{\partial x^{\frac{1}{3}} }. $ 
For $w_{t,\frac{1}{3}}=\frac{1}{3}$ by using remark \eqref{BB} we can write
$$u_{t}+2\sqrt[3]{xu_{x}}+u=x^{2} \Longrightarrow u_{t}+2\sqrt[3]{x} w_{t,\alpha} x^{1-\alpha}u_{x}+u=x^{2} \Longrightarrow 
u_{t}+\frac{2}{3}xu_{x}+u=x^{2}.
$$
We solve this equation by taking Laplace transform of equation respect to $ t $, we denote by $U(x,s)$ the Laplace of $ u(x,t) $ respect to $ t $, we have the following equation   
\begin{equation}\label{tropsol2}
\begin{cases}
sU(x,s)-U(x,0)+ \frac{2}{3}x U_{x}(x,s)+U(x,s)={x^2}{s}\\
U(t,0)=0,\\
U(0,x)=0.
\end{cases}
\end{equation}
then 
$$ U_{x}+\frac{3+3s}{2x}U=\frac{3x}{2s}\Longrightarrow U(x,s)=\frac{3x^2}{s(3s+7)}+c(s)x^{\frac{-3-3s}{2}}.$$
By substituting $U(0,x)=0$ we have $ c(s)=0 $, then $ U(x,s)=x^2(\frac{3}{7s}+\frac{9}{7(3s+7)})$. The solution of equation is $u(x,t)=\frac{x^2}{7}(1-e^{\frac{-7}{3}}) .$  
\end{example}
\begin{example}
We consider the partial fractional differential equation of second order; 
\begin{equation}\label{22}
 \sqrt[5]{u_{xt}}+2\frac{u}{x}=0, 
 \end{equation}
where $ u(x,t) $ be a $ \frac{1}{5}- $fractional partial differentiable function of second order respect to $ t,x $. 
For $ w_{x,\frac{1}{5}}=x^2, w_{t,,\frac{1}{5}}=\frac{1}{\sqrt[3]{t}} $ by using remark \eqref{CC} we have
\begin{equation}\label{11}
x^2 \sqrt[5]{x^4} \sqrt[5]{t^4} u_{xt}+2\sqrt[3]{t}u=0.
\end{equation}
We consider a solution of this differential equation with the form $ u(x,t)=f(x)g(t) $ such that $ f $ a function depends on $ x $ and $ g $ a function depends on
$ t $. By substituting $ u(x,t) $ in the equation \eqref{11} we have
\begin{equation}\label{44}
x^2 \sqrt[5]{x^4} \sqrt[5]{t^4} Df.Dg+3\sqrt[3]{t} fg=0. 
\end{equation}
We can write the equation \eqref{44} a form that divide the functions of $ t $ and $ x $:
\begin{equation}\label{33}
\frac{x^2 \sqrt[5]{x^4} Df}{f}=-\frac{2\sqrt[3]{t} g}{\sqrt[5]{t^4}Dg},
\end{equation}
 two sides of the equality \eqref{33} is a constant $ k $. We have 
\begin{equation}\label{tropsol2}
 \begin{cases}
\frac{Df}{f}=\frac{k}{\sqrt[5]{x^{14}}}\rightarrow   Ln f= k\sqrt[15]{x}+c_{1} \rightarrow f=\exp( k\sqrt[15]{x}+c_{1})\\
\frac{Dg}{g}=\frac{-2}{k\sqrt[15]{t^{7}}}\rightarrow Ln g= \frac{-2\sqrt[15]{t^8}}{k}+c_{2}\rightarrow g= \exp(\frac{-2\sqrt[15]{t^8}}{k}+c_{2})
\end{cases}
\end{equation}
The solution of the equation \eqref{22} has the form 
$$ u(x,t)= \exp( k\sqrt[15]{x}+c_{1})\exp(\frac{-2\sqrt[15]{t^8}}{k}+c_{2})=c\exp(k\sqrt[15]{x}+\frac{-2\sqrt[15]{t^8}}{k}).$$
\end{example}

\section{Conclusion}
We define a generalized fractional derivative(GFD). We show that the previous derivatives
are particular cases. We show how it is possible to have infinite
fractional derivatives
with their algebra. We present the fractional differential ring, the
fractional partial
derivatives and their applications.
\section*{Acknowledgement}
This paper was supported by PRODEP of Mexico for Zeinab Toghani postdoctoral position.

\newpage
\bibliographystyle{amsalpha} \addcontentsline{toc}{chapter}{Bibliografia}
\bibliography{bibliografia_zeinab}

\end{document}